\documentstyle[12pt]{article}
\newcommand{\supp}{\mbox{supp}}
\newcommand{\Span}{\mbox{span}}

\def \N{\hbox{$I\hskip -4pt N$}} \def\sN{\hbox{$\sc I\hskip -3pt N$}}
\def \Z {\hbox{$Z\hskip -5.2pt Z$}}
\def\sZ{\hbox{$\sc Z\hskip -4.2pt Z$}}
\def \Q{\hbox{$Q\hskip -6pt \vrule height 6pt depth 0pt\hskip 6pt$}}

\def \C{\hbox{$C\hskip -7pt \vrule height 6pt depth 0pt \hskip 6pt$}}
\def \sC{{\hbox{$\sc C\hskip -5pt \vrule height 5pt depth 0pt \hskip 6pt$}}}
\def\qed{\hfill \hfill \ifhmode\unskip\nobreak\fi\ifmmode\ifinner
         \else\hskip5pt\fi\fi
 \hbox{\hskip5pt\vrule width4pt height6pt depth1.5pt\hskip 1 pt}}
\def\a{\alpha}
\def\b{\beta}

\def\l{\lambda}

\def\j{\setminus}

\def\Vir{\hbox{\rm Vir}}
\def\rank{\hbox{\rm rank}}

\def\sc{\scriptstyle}

\def\cl{\centerline}

\begin{document}
\par\
\par\
\par

\cl{\large{\bf  Classification of irreducible Harish-Chandra
modules over}}

\cl{\large{\bf  generalized Virasoro algebras}\footnote{
\indent \hskip .3cm  Research supported by NSERC, and the NSF  of
China  (Grants 10371120  and 10431040).}}

\par
\vskip 15pt
 \centerline{Xiangqian
Guo, Rencai Lu and Kaiming Zhao}

\par
\vskip 15pt

\cl{Institute of Mathematics}

\cl{Academy of Mathematics and System Sciences}

\cl{Chinese Academy of Sciences} \cl{Beijing 100080, P. R. China}

\vskip 10pt \centerline{ and } \vskip 10pt
 \centerline{
Department of Mathematics}

\centerline{  Wilfrid Laurier University}

\centerline{  Waterloo, ON, Canada N2L 3C5}

\vskip 5pt \centerline{  Email: kzhao@wlu.ca}

%\centerline{  Email: kzhao@math08.math.ac.cn}
\par\vskip 10pt
%\vs{5pt}
\begin{abstract}
{Let $G$ be an arbitrary additive subgroup of $\C$ and $\Vir[G]$
the corresponding generalized Virasoro algebra. In the present
paper, irreducible weight modules with finite dimensional weight
spaces over $\Vir[G]$ are completely determined. The
classification strongly depends on the index group $G$. If $G$
does not have a direct summand $\Z$, then  such irreducible
modules over $\Vir[G]$ are only modules of intermediate series
whose weight spaces are all $1$-dimensional. Otherwise, there is
one more class of modules  which  are constructed by using
intermediate series modules over a generalized Virasoro subalgebra
$\Vir[G_0]$ of $\Vir[G]$ for a direct summand $G_0$ of $G$ with
corank $1$.}
\end{abstract}
\vskip 30pt

\vskip .5cm

\break
\par
\cl{{\bf \S1. Introduction}}
\par
\vskip .3cm

The Virasoro algebra theory has been widely used in many physics
areas and other mathematical branches, for example, string theory,
2-dimensional conformal field theory, differential geometry,
combinatorics, Kac-Moody algebras, vertex algebras, and so on.
\vskip 5pt

The generalized Virasoro algebras were first introduced and
studied by mathematicians and mathematical physicists J. Patera
and H. Zassenhaus [PZ] in 1991. Because of  its own interest and
the close relation between the theory of generalized Virasoro
algebras and physics, this theory has attracted extensive
attentions of mathematicians and physicists, particularly, the
representation theory of generalized Virasoro algebras has been
developed rapidly in the last decade. Let us first recall the
definitions of these Lie algebras. \vskip 5pt

In this paper we denote by $\C$,  $\Z$, $\Z_+$ and $ \N$ the set
of  complex numbers, integers, nonnegative integers and positive
integers respectively.

\vskip 5pt The {\bf Virasoro algebra $\Vir:=\Vir[\Z]$} (over $\C$)
is  the Lie algebra with the basis
$\bigl\{C,d_{i}\bigm|i\in\Z\bigr\}$ and the Lie brackets defined
by

$[d_m,d_n]=(n-m)d_{m+n}+\delta_{m,-n}\frac{m^{3}-m}{12}C,\qquad\forall
m,n\in \Z,$

$[d_m,C]=0, \qquad\forall m\in \Z.$ \vskip 5pt

The structure theory of Harish-Chandra modules over the Virasoro
algebra is  developed fairly well. For details, we refer the
readers to [Ma], [MP], [MZ], the book [KR] and the references
therein. In particular, the classification of irreducible
Harish-Chandra modules was obtained in [Ma], while indecomposable
modules was studied in [MP]. This classification was recently used
to give the classification of irreducible weight modules over the
twisted Heisenberg-Virasoro algebra [LZ2]. \vskip 5pt

Patera and Zassenhaus [PZ] introduced the {\bf generalized
Virasoro algebra} $\Vir[G]$ for any additive subgroup $G$ of $\C$
from the context of mathematics and physics. This Lie algebra can
be obtained from $\Vir$ by replacing the index group $\Z$ with $G$
(see Definition 2.1). This Lie algebra $\Vir[G]$ is called a {\bf
rank n Virasoro algebra} (or a {\bf higher rank Virasoro algebra}
if $n\geq 2$) if $G\simeq \Z^n$. \vskip 5pt

Representation theory of generalized Virasoro algebras have been
extensively studied in recent years. Mazorchuk [M] proved that all
irreducible Harish-Chandra modules over $\Vir[\Q]$ are
intermediate series modules (where $\Q$ is the field of rational
numbers). In [HWZ], a criterion for the irreducibility of Verma
modules over the generalized Virasoro algebra $\Vir[G]$ was
obtained. In [S1], [S2] and [SZ], the irreducible Harish-Chandra
modules over the generalized Virasoro algebras were investigated.
[BZ] constructed a new class of irreducible Harish-Chandra modules
over some generalized Virasoro algebras. Recently, a complete
classification of irreducible Harish-Chandra modules over higher
rank Virasoro algebras was given by  the last two authors of the
present paper [LZ1]. \vskip 5pt

In this paper, we give the classification of irreducible
Harish-Chandra modules over any generalized Virasoro algebra.

%{\bf Theorem 1.1.}\ {\it Let $G$ be any additive subgroup of $\C$.
%Then any nontrivial irreducible Harish-Chandra module over
%$Vir[G]$ is either of intermediate series or isomorphic to some
%$V(\a,\b,G_0,b)$ for some $\a,\b\in \C,\, b\in G$ and a subgroup
%$G_0$ of $G$ such that $G=G_0\oplus \Z b$.}
\vskip 5pt

The paper is organized as follows. In Section 2, for the reader's
convenience, we collect some results and give some notations from
[M], [MP], [SZ] and [LZ1] for later use. In Section 3,  we widely
use the results in [M], [LZ1] and [SZ] to complete the
classification of irreducible weight modules with finite
dimensional weight spaces over $\Vir[G]$. The classification
strongly depends on the index group $G$. If $G$ does not have a
direct summand $\Z$, then  such irreducible modules over $\Vir[G]$
are only modules of intermediate series whose weight spaces are
all $1$-dimensional. Otherwise, there is one more class of modules
which are constructed by using intermediate series modules over a
generalized Virasoro subalgebra $\Vir[G_0]$ of $\Vir[G]$ for a
direct summand $G_0$ of $G$ with corank $1$ (Theorem 3.8). We
break the proof into seven lemmas.

\vskip 5pt Throughout this paper, a subgroup always means an
additive subgroup if not specified. For any $a\in \C$ and
$S\subset \C$, we denote $a+S=\{a+x| x\in S\}$ and $aS=Sa=\{ax|
x\in S\}$.

\vskip .5cm

\par
\cl{{\bf \S2. Modules over generalized Virasoro algebras }}
\par
\vskip .2cm

First we give the definition of the generalized Virosoro algebras.

{\bf Definition 2.1.} {\it Let $G$ be a nonzero additive subgroup
of $\C$. The {\bf generalized Virasoro algebra $\Vir[G]$} (over
$\C$) is the Lie algebra with the
 basis $\bigl\{C,d_{x}\bigm|x\in G\bigr\}$ and
the Lie bracket defined by $$ [d_{x},
d_y]=(y-x)d_{x+y}+\delta_{x,-y}\frac{x^{3}-x}{12}C,\qquad
\forall\,\,x,y\in G,$$
$$[C, d_{x}]=0,\,\,\forall\,\,x \in G.$$ }

It is clear that $\Vir[G]\cong \Vir[aG]$ for any $a\in \C^*=\C\j
\{0\}$. Then for any $x\in G^*:= G\j\{0\}$, $\Vir[x\Z]$ is a Lie
subalgebra of $\Vir[G]$ isomorphic to $\Vir=\Vir[\Z]$, the
classical Virasoro algebra.

\medskip

A $\Vir[G]$-module $V$ is called {\bf trivial} if ${\Vir}[G]V=0$.
For any $\Vir[G]$-module $V$ and $c,\lambda\in \C$,
$V_{\lambda,c}:=\bigl\{v\in V\bigm|d_{0}v=\lambda v, Cv=cv\bigr\}$
is called the {\bf weight space} of $V$ corresponding to the
weight $(\lambda,c)$. When $C$ acts as a scalar $c$ on the whole
module $V$, we shall simply write $V_{\lambda}$ instead of
$V_{\lambda,c}$. In the rest of this section, all modules
considered are such modules.

\medskip
A $\Vir[G]$-module $V$ is called a {\bf weight module} if $V$ is
the sum of its weight spaces, and a weight module is called a {\bf
Harish-Chandra module} if all the weight spaces are finite
dimensional. For a weight module $V$, we define $\supp
V:=\bigl\{\lambda\in \C\bigm|V_{\lambda}\neq 0\bigr\}$, which is
generally called the {\bf weight set} (or the {\bf support}) of
$V$. Given a weight module $V$ and any subset $S\subset \C$, we
denote $V_{S}=\bigoplus_{x\in S} V_{x}$, where $V_{x}=0$ for
$x\not\in \supp V$.

\medskip
Let $V$ be a module and $W'\subset W$ are submodules of $V$. The
module $W/W'$ is called a {\bf sub-quotient} of $V$. If $W'=0$ we
consider that $W=W/W'$.

\medskip
 Let $V$ be a weight module over $\Vir[G]$. $V$ is said to be
{\bf uniformly bounded}, if there exists $N\in \N$ such that $\dim
V_{x}<N$ for all $x\in \supp V$.

\medskip
Fix a total order "$\succeq$" on $G$ which is compatible with the
addition, i.e., $x\succeq y$ implies $x+z\succeq y+z$ for any
$x,y,z\in G$. Let
$$ G^{+}:=\bigl\{x\in G\bigm|x\succ 0\bigr\},\quad
G^{-}:=\bigl\{x\in G\bigm|x\prec 0\bigr\},
$$
$$ \Vir[G]^{+}:=\sum_{x\in G^+}\C d_x,\quad
\Vir[G]^{-}:=\sum_{x\in G^-}\C d_x.
$$
Let $V$ be a weight module over $\Vir[G]$. A vector $v\in V_{\l,
c},\, \l\in \supp V$, $c\in\C$, is called a {\bf highest weight
(resp. lowest weight) vector} if $\Vir[G]^{+} v=0$ (resp.
$\Vir[G]^{-} v=0$). $V$ is called a {\bf highest weight (resp.
lowest weight) module} with highest weight (resp. lowest weight)
$(\l,c)$ if there exists a nonzero highest (lowest, resp.) weight
vector $v\in V_{\lambda,c}$ such that $V$ is generated by $v$. For
the natural total order on $\Z$, the irreducible highest weight
$\Vir[\Z]$-module with highest weight $(\l,c)$ is generally
denoted by $V(c,\l)$.
\medskip

Now we give another class of weight modules over $\Vir[G]$, i.e.,
the {\bf modules of intermediate series $V(\a,\b,G)$}. For any
$\a,\b\in \C$, the module $V(\a,\b,G)$ has a basis $\{v_x\,|\,x\in
G\}$ with actions of $\Vir[G]$ given by:

$$ C v_y=0,\,\,\,\, d_x v_y=(\a+y+x
\b)v_{x+y},\,\,\forall\,\,x,y\in G.$$

One knows from [SZ] that $ V(\a,\b,G)$ is reducible if and only if
$\a\in G$ and $\b\in\{0,1\}$. By $V'(\a,\b,G)$ we denote the
unique nontrivial irreducible sub-quotient of $V(\a,\b,G)$. Then
$\supp(V'(\a,\b,G))=\a+G$ or $\supp(V'(\a,\b,G))=G\j\{0\}$. We
also refer $V'(\a,\b,G)$ as {\bf intermediate series modules}. The
following result is due to Su and Zhao.

\medskip {\bf Theorem 2.2 ([SZ, Theorem 4.6]).} {\it Let $V$ be a
nontrivial irreducible Harish-Chandra module over $\Vir[G]$ with
all weight spaces $1$-dimensional. Then $V \cong V'(\a,\b,G)$ for
some $\a,\b\in\C$.}

\medskip
This result for the classical Virasoro algebra is due to Kaplansky
[K1]. The following classification of irreducible Harish-Chandra
modules over the classical Virasoro algebra was obtained by
Mathieu.
\medskip

{\bf Theorem 2.3 ([Ma]).}\ {\it Every irreducible Harish-Chandra
module over $\Vir$ is either a highest weight module, a lowest
weight module, or a module of intermediate series.}

\medskip
We say that a $\Vir[\Z b]$-module $W$ is {\bf positively
truncated} ({\bf negatively truncated}) relative to $b$ if for any
$\l\in \supp V$, there exists some $x_0\in \Z$ such that
$\supp(V_{\l+\sZ b})\subset \{\l+xb|x\leq x_0\}$ (resp.
$\supp(V_{\l+\sZ b})\subset \{\l+xb|x\geq x_0\}$). The following
result from Martin and  Piard will be useful to our later proofs.

\medskip {\bf Theorem 2.4 ([MP]).}\  {\it Every Harish-Chandra
$\Vir$-module $V[\Z b],\, b\in \C$ which has neither trivial
submodules nor trivial quotient modules can be decomposed as a
direct sum of three submodules $V=V^+\oplus V^0\oplus V^-$, where}
$V^+$ {\it is positively truncated,} $V^-$ {\it is negatively
truncated and $V^0$ is uniformly bounded.}

\medskip
Now we come back to the generalized Virasoro algebras. Marzorchuck
proved

\medskip {\bf Theorem 2.5 ([M]).} {\it Any nontrivial irreducible
Harish-Chandra module over  $\Vir[\Q]$ is a module of intermediate
series.}

\medskip In fact, he has proved the following theorem as he
remarked at the end of his paper [M]:

\medskip {\bf Theorem 2.5$'$.} {\it Let $G$ be an infinitely
generated additive subgroup of $a\Q$ for some $a\in \C$. Then any
nontrivial irreducible Harish-Chandra module over $\Vir[G]$ is a
module of intermediate series.}

\medskip
Now we assume that $G= G_0\oplus \Z b\subset\C$ where $0\neq b \in
\C$ and $G_0$ is a nonzero subgroup of $\C$. (Note that some $G$
does not possess this property, for example, $\Q$). Set
%\medskip
$$\Vir[G]_+=\bigoplus_{x\in G_0,k\in \sZ^+}\C d_{x+kb}\oplus\C
C$$ and $$\Vir[G]_{++}=\bigoplus_{x\in G_0,k\in \sN}\C d_{x+kb}.$$
%\medskip
Given any $\a,\b\in\C$, let $V'(\a,\b,G_0)$ be the module of
intermediate series over $\Vir[G_0]$. We extend the
$\Vir[G_0]$-module structure on $V'(\a,\b,G_0)$ to a
$\Vir[G]_+$-module structure by defining
$\Vir[G]_{++}V'(\a,\b,G_0)=0$. Then we obtain the induced
$\Vir[G]$-module
$$ M(b,G_0,V'(\a,\b,G_0))%=\Ind_{\Vir[G]_+}^{\Vir[G]}V'(\a,\b,G_0)
=U(\Vir[G])\otimes_{U(\Vir[G]_+)}V'(\a,\b,G_0),$$
 where $U(\Vir[G])$ and $U(\Vir[G]_+)$ are the universal
enveloping algebras of $\Vir[G]$ and $\Vir[G]_+$, respectively.

\medskip The $\Vir[G]$-module $M(b,G_0,V'(\a,\b,G_0))$ has a unique
maximal proper submodule $J$. Then we obtain the irreducible
quotient module $$V(\a,\b, b,G_0):=M(b,G_0,V'(\a,\b,G_0))/J.$$ It
is clear that this module is uniquely determined by $\a,\b, b$ and
$G_0$ and that ([Lemma 3.8, LZ1])

\medskip \cl{$\supp V(\a,\b,b,G_0)=\ \a-\Z^+b+G_0$\,\,\,\ or
$\,\,\,(-\Z^+b+G_0)\j\{0\}.$} \medskip

It was proved that $V(\a,\b,b,G_0)$ is a Harish-Chandra module
over $\Vir[G]$:

\medskip {\bf Theorem 2.6 ([BZ, Theorem 3.1]).} {\it The $\Vir[G]$-modules
$V(\a,\b,b,G_0)$ are Harish-Chandra modules. More
precisely, $\dim V_{-i b+\a+x}\le (2i+1)!!$ for all $ i\in \N,\,
x\in G_0.$}

\medskip
 From this theorem, we
easily deduce the following corollary which will be used
frequently in our later proofs.

\medskip {\bf Corollary 2.7.} {\it Let $V=V(\a,\b, b,G_0)$ and $i\in \Z^+$.
Then for any subgroup $G'$ of $G$, the $\Vir[G']$ module
$V_{\a-ib+G'}$ is uniformly bounded if and only if $G'\subset
G_0$.}

\medskip
 The classification of irreducible Harish-Chandra modules over higher rank Virasoro
algebra was obtained by the last two author of the present paper.

\medskip {\bf Theorem 2.8 ([LZ1]).}  {\it Let $G$ be an additive
subgroup of $\C$ such that $G\cong \Z^n$ for some $n\in \N$ with
$n>1$. Then any nontrivial irreducible Harish-Chandra module over
$\Vir[G]$ is either of intermediate series or isomorphic to some
$V(\a,\b,G_0,b)$ for some $\a,\b\in \C,\, b\in G$ and a subgroup
$G_0$ of $G$ such that $G=G_0\oplus \Z b$.}

%\medskip We note here the fact that, for any finite subset $S $ of
%$\C$, the subgroup $\langle S\rangle$ of $\C$ generated by $S$ is
%isomorphic to $\Z^n$ for some $n\in\N$. We will use this result
%frequently in the sequel.

\vskip .5cm
\par
\cl{{\bf \S3. Classification of irreducible Harish-Chandra
modules}}

\cl{{\bf over generalized Virasoro algebras}}
\par
\vskip .2cm

In this section we give a classification of the irreducible
Harish-Chandra modules over the generalized Virasoro algebras. Let
us proceed with the convention that $U(G)=U(\Vir[G])$, the
enveloping algebra, and that $U(G)_a=\{ y\in U(\Vir[G])|\,
[d_0,y]=ay\}, a\in \C$, for any additive subgroup $G$ of $\C$.

\medskip
We recall the concept: the rank of a subgroup $A$ of $\C$ from
[K2]. The {\bf rank} of $A$, denoted by rank$(A)$, is the maximal
number $r$ with $g_1,\cdots,g_r\in A\j\{0\}$ such that $\Z
g_1+...+\Z g_r$ is a direct sum. If such an $r$ does not exist, we
define $\rank(A)=\infty$.

\medskip From now on we fix a subgroup $G$
of $\C$, and a
 nontrivial irreducible Harish-Chandra module $V$ over
$\Vir[G]$. Because of the classifications in Theorems 2.2 and
2.5$'$, we may also {\bf assume that rank$G\ge 2$.}

%$U(G)_0=U(\Vir[G])_0=\bigl\{y\in U(\Vir[G])\bigm| [d_0,y]=0\bigr\}$.

\medskip
{\bf Lemma 3.1.}  {\it For any finite subset $I$ of $\supp V$,
there is a subgroup $G_I$ of $G$ such that}

(a) {\it $G_I\cong \Z^k$ for some $k\in \N$,}

(b) {\it $U(G_I)V_{\mu}=U(G_I)V_{\mu'}$ and $\mu-\mu'\in G_I$ for
any $\mu, \mu'\in I$,}

(c) {\it $V_{\mu}$ is an irreducible $U(G_I)_0$ module for any
$\mu\in I$.}

\medskip

{\it Proof.} Assume that $I=\{\mu_i| i=1,2,...,s \}$. Since $V$ is
an irreducible $U(G)$-module and $V_{\mu_i}$ are finite
dimensional, then for any $\mu_i,\mu_j\in I$, there are elements
$y_{i,j}^{(1)},...,y_{i,j}^{(d_{ij})}\in U(G)_{\mu_i-\mu_j}$ such
that $V_{\mu_i}= \sum_{t=1}^{d_{ij}}y_{i,j}^{(t)}V_{\mu_j}$, where
$d_{ij}\in\N$.

 For any $\mu_i\in \supp V$, since $V$ is irreducible we see that  $V_{\mu_i}$ is an irreducible
$U(G)_0$-module. Let $\phi_{\mu_i}:\ U(G)_0\longrightarrow
gl(V_{\mu_i})$ be the representation of $U(G)_0$ in $V_{\mu_i}$,
where $gl(V_{\mu_i})$ is the general linear Lie algebra associated
with the vector space $V_{\mu_i}$. Since $V_{\mu_i}$ and hence
$gl(V_{\mu_i})$ is finite dimensional, then $\phi_{\mu_i}(U(G)_0)$
is finite dimensional. Thus there exist $y_{i1},...,y_{im_i}\in
U(G)_0$ such that
$\Span_{\sC}\{\phi_{\mu_i}(y_{i1}),...,\phi_{\mu_i}(y_{im_i})\}=\phi_{\mu_i}(U(G)_0)$.

For the finitely  many elements: $y_{i1},...,y_{im_i},\,\,\,
y_{i,j}^{(1)},...,y_{i,j}^{(d_{ij})};$ $1\le i,j\le s$, thanks to
the PBW Theorem there are finitely  many elements $g_1,...,g_n\in
G$ such that $y_{i,j} ,y_{i,j}^{(k)}\in U(G_I)$, where $G_I$ is
the subgroup of $G$ generated by $g_1,...,g_n$.

Since $G_I$ is a finitely generated torsion free abelian group, it
is clear that $G_I\cong\Z^k$ for some $k\in \N$, (a) follows. By
the construction of $G_I$, we know that $V_{\mu_i}\subset
U(G_I)V_{\mu_j}$. Hence $U(G_I)V_{\mu_i}=U(G_I)V_{\mu_j}$ for any
$\mu_i,\mu_j\in I$. It is also clear that $\mu-\mu'\in G_I$ for
any $\mu, \mu'\in I$. Thus (b) follows.

We now prove that $V_{\mu_i}$ is an irreducible $U(G_I)_0$-module
for any $\mu_i\in I$. Suppose that $N$ is a nontrivial
$U(G_I)_0$-submodule of $V_{\mu_i}$. For any element $y\in
U(G)_0$, there are some $a_j\in \C$   such that
$\phi_{\mu_i}(y)=\sum_{j=1}^{m_i}a_j\phi_{\mu_i}(y_{ij})$. Then
$yN=\sum_{j=1}^{m_i}a_j\phi_{\mu_i}(y_{ij})N=(\sum_{j=1}^{m_i}a_jy_{ij})N\subset
N$. That is, $N$ is a nonzero $U(G)_0$-submodule of $V_{\mu_i}$,
forcing $N=V_{\mu_i}$. Hence, $V_{\mu_i}$ is an irreducible
$U_0(G_I)$-module for any $\mu_i\in I$, and (c) follows. \qed

\medskip
{\bf Lemma 3.2.}  {\it Let $I$ and $G_I$ be the same as in Lemma
3.1 and let $G'$ be a subgroup of $G$ that contains $G_I$. Then
for any $\l\in I$, $V_{\l+G'}$ has a unique irreducible
$\Vir[G']$-sub-quotient $V'$ with $\dim V'_{\mu}=\dim
V_{\mu},\,\forall \mu\in I$.}

\medskip
{\it Proof.}\ Since $G_I\subset G'$,  from (b) and (c) of Lemma
3.1, we know that,

(b$'$) $U(G')V_{\mu}=U(G')V_{\mu'}$ for any $\mu, \mu'\in I$,
which we denote by $W$.

(c$'$) $V_{\mu}$ is an irreducible $U_0(G')$ module for any
$\mu\in I$,

\noindent Clearly, $W\subset V_{\l+G'},\,\, \forall \l\in I$.

Suppose that  $W'$ is a proper $\Vir[G']$-submodule of $W$. For
any $\mu\in I$, it is clear that $W'_{\mu}$ is a proper $U(G')_0$
submodule of $W_{\mu}$, forcing $W'_{\mu}=0$, that is, $W'$
trivially intersects $W_{\mu}$  for any $\mu\in I$. Now let
$V'=U(G')X/Y$, where $X=\oplus_{\mu\in I} V_{\mu}$ and $Y$ is the
sum of all $U(G')$-submodules $Y'$ of $U(G')X$ such that
$Y'_{\mu}=0$ for all $\mu\in I$. Then $V'$ is as desired.\qed

\medskip
{\bf Theorem 3.3.} (a) {\it If $V$ is uniformly bounded, then $V$
is of intermediate series.}

(b) {\it A Harish-Chandra module $W$ over $\Vir[G]$ with
$\supp(W)\subset \l+G$ for some $\l\in\C$ is uniformly bounded if
and only if $\dim W_{\l}$ $=\dim W_{\mu}$, for all $\l,\mu\in
\supp V\j\{0\}$.}

\medskip
{\bf Remark.} This theorem holds also for any rank $1$ group $G$.

\medskip
{\it Proof.} (b) follows directly from (a). To prove (a), it
suffices to show that $\dim V_{\l}=1$, for all $\l\in \supp V$, by
Theorem 2.2.

Now suppose that $\dim V_{\l}\geq 2$ for some $\l\in \supp V$.
Using Lemma 3.1 for $I=\{\l\}$, we have the subgroup $G_I$ of $G$
described there. Then $V_{\l}$ is an irreducible
$U(G_I)_0$-module. Consider the uniformly bounded
$Vir[G_I]$-module $V_{\l+G_I}$, then it has an irreducible
uniformly bounded sub-quotient   $V'$ with $\dim V'_{\l}=\dim
V_{\l}\ge2$ since $V_{\l}$ is an irreducible $U(G_I)_0$-module. by
Theorem 2.8, $\dim V'_{\l}$ should be not larger  than  $1$, a
contradiction. Thus we have that $dim V_{\l}=1,\,\, \forall \l\in
\supp V$. The proof is completed.\qed

\medskip
{\bf Lemma 3.4.}  {\it For any} $\l\in \supp V$ {\it and} $g\in
G\setminus\{0\}$, {\it if the $Vir[\Z g]$-module $V_{\l+\sZ g}$
has a nontrivial uniformly bounded $\Vir[\Z g]$-sub-quotient, then
$V_{\l+\sZ g}$ itself is uniformly bounded.}

\medskip
{\it Proof.} To the contrary, suppose that $V_{\l+\sZ g}$ is not
uniformly bounded. Then there exist some $\mu_1,\mu_2\in \supp
(V_{\l+\sZ g})\j\{0\}$ with $\dim V_{\mu_1}\neq \dim V_{\mu_2}$.

Applying  Lemma 3.1 to $I=\{\mu_1,\mu_2\}$, we have the subgroup
$G_I$ of $G$ described there, and furthermore  we may assume that
$\rank G_I>1$. Let $G'=G_I+\Z g$. Then by Lemma 3.2, the
$\Vir[G']$-module $V_{\l+G'}$  has a unique irreducible
sub-quotient $V'$ with $\dim V'_{\mu_i}=\dim V_{\mu_i}, i=1,2$.

Since $\dim V'_{\mu_1}\neq \dim V'_{\mu_2}$, Theorem 3.3 ensures
that $V'$ is not uniformly bounded. By Theorem 2.8, $V'$ must be
of the form $V(\a,\b,G'_0,b)$, for some $\a,\b\in \C,\, b\in G'$
and a subgroup $G'_0$ of $G'$ with $G'=G'_0\oplus \Z b$. Since
$\dim V'_{\mu_1}\neq \dim V'_{\mu_2}$ and $\mu_1\ne0\ne\mu_2$,
$V'_{\l+\sZ g}$ is not uniformly bounded by Theorem 3.3. Then by
Corollary 2.7, we know that $g\notin G'_0$. Thus the $\Vir[\Z
g]$-module $V'_{\l+\sZ g}$ is positively truncated relative to
$g$.

From the definition of $V'$ we know that $V'_{\l+\sZ g}$ is a
$\Vir[\Z g]$-sub-quotient of $V_{\l+\sZ g}$, \break say,
$V'_{\l+\sZ g}=W/W'$ where $W'\subset W$ are $\Vir[\Z
g]$-submodules of $V_{\l+\sZ g}$. Since $\dim V_{\mu_1}=\dim
V'_{\mu_1}=\dim W_{\mu_1}$, then $W'_{\mu_1}=0$ and $(V_{\l+\sZ
g}/W)_{\mu_1}=0$. Note that $\mu_1\ne0$. Therefore, $V_{\l+\sZ
g}/W$, $W/W'$, $W'$ all do not have uniformly bounded $\Vir[\Z
g]$-sub-quotient. Thus $V_{\l+\sZ g}$ does not have any nontrivial
uniformly bounded $\Vir[\Z g]$-sub-quotient, a contradiction.\qed

%Then we have that for any $v\in V'_{\mu_i}=V_{\mu_i},\, i=1,2$,
%$d_{kg}v=0$ or $d_{-kg}v=0$ for sufficiently large $k\in \N$. {\bf
%This is not right!} That is, $V_{\l+\sZ g}$ does not contain any
%nontrivial uniformly bounded $\Vir[\Z g]$-sub-quotient, a
%contradiction.\qed

\medskip
{\bf Lemma 3.5.} {\it Assume that} $\l\in\supp V\j\{0\}$ {\it and
that $G_1, G_2$ are any subgroups of $G$. If both $V_{\l+G_1}$ and
$V_{\l+G_2}$ are uniformly bounded (as $\Vir[G_1]$-module and
$\Vir[G_2]$-modules, resp.), then $V_{\l+G_1+G_2}$ is a uniformly
bounded $\Vir[G_1+G_2]$-module.}
\medskip

{\it Proof.}\ Thanks to Theorem 3.3, we may assume that $\dim
V_{\mu}=m$ for nonzero $ \mu\in( \l+G_1)\cup (  \l+G_2)$.
%Denote $m=\max\{m,\dim V_0 \}$.

To the contrary, we suppose that there are some $g_i\in G_i,
i=1,2$ such that $\dim V_{\l+g_1+g_2}\neq m$ with
$(\l+g_1)(\l+g_2)(\l+g_1+g_2)\neq 0$. Applying Lemma 3.1 to
$I=\{\l+g_1, \l+g_2, \l+g_1+g_2\}$, we have the subgroup $G_I$ as
described in Lemma 3.1, and furthermore  we may assume that $\rank
G_I>1$. Note that $  g_1, g_2\in G_I$. Lemma 3.2 ensures that the
$\Vir[G_I]$-module $V_{\l+G_I}$ has a unique irreducible
sub-quotient $V'$ such that $\dim V'_{\mu}=\dim V_{\mu}$ for all
$\mu\in I$. In particular, $\dim V'_{\l+g_1}\ne \dim
V_{\l+g_1+g_2}$. Noting that $(\l+g_1)(\l+g_1+g_2)\ne0$, by
Theorem 3.3 we know that $V'$ is not uniformly bounded.
 By Theorem 2.8, $V'\cong V(\a,\b,G_0,b)$ for
some $\a,\b\in \C,\, b\in G_I$ and a subgroup $G_0$ of $G_I$ with
$G_I=G_0\oplus \Z b$. The fact that $V'_{\l+\sZ g_1}$ and
$V'_{\l+\sZ g_2}$ are both uniformly bounded implies $g_1, g_2\in
G_0$. Thus $g_1, g_2,g_1+g_2\in G_0$. By Corollary 2.7 , we deduce
that $V'_{\l+\sZ g_1+\sZ g_2}$ is uniformly bounded. Then by
Theorem 3.3, we see that $\dim V_{\l+g_1+g_2}=\dim V_{\l+g_1}=m$,
a contradiction. Hence, $\dim V_{\l+x}=m$ for any $x\in G_1+G_2$
with $\l+x\ne0$. Therefore $V_{G_1+G_2}$ is uniformly bounded.\qed

%and hence $\Z g_1+\Z g_2\subset G'_0$. Then we have that $V'_{\l+g_1+\Z (g_1+g_2)}$ is
%uniformly bounded, which is a $\Vir[\Z(g_1+g_2)]$-sub-quotient of $V[\l+g_1+\Z(g_1+g_2)]$, and hence
%that $V[\l+\Z(g_1+g_2)]$ is uniformly bounded by Lemma 3.4, contradiction .

\medskip
{\bf Lemma 3.6.}  {\it For any} $\mu\in \supp V\j\{0\}$, {\it
there exists a unique maximal subgroup $G_{\mu}$ of $G$ such that
$V_{\mu+G_{\mu}}$ is a uniformly bounded $\Vir[G_{\mu}]$-module.
Furthermore,}

  (a) {\it $G_{\mu_1}=G_{\mu_2}$ for any}  $\mu_1,\mu_2\in \supp
V\j\{0\}$, {\it which we denote by $G^{(0)}$,}

  (b) {\it either $G^{(0)}=G$ or $G\cong G^{(0)}\oplus \Z b$ for some
$b\in G$.}

\medskip
{\it Proof.}\ \ The existence and uniqueness of $G_{\mu}$ for
$\mu\in \supp V\j\{0\}$ follows easily  from Lemma 3.5.

\medskip
(a)\ Fix  $\mu_1\neq\mu_2\in \supp V\j\{0\}$ and  $g_1\in
G_{\mu_1}\j\{0\}$. Let $I=\{\mu_1,\mu_2\}$ and take $G_I$ the same
subgroup    described as in lemma 3.1, and furthermore  we may
assume that $\rank G_I>1$. Set $G'=G_I+\Z (\mu_1-\mu_2)+\Z g_1$.
By Lemma 3.2, the $\Vir[G']$-module $V_{\mu_2+G'}$ has a unique
irreducible sub-quotient $V'$ with $\dim V'_{\mu_i}=\dim
V_{\mu_i}, i=1,2$. Clearly, $V'$ is nontrivial.

If $V'$ is uniformly bounded, then $V'_{\mu_2+\sZ g_1}$ is a
uniformly bounded $\Vir[\Z g_1]$-sub-quotient of $V_{\mu_2+\sZ
g_1}$.

If $V'$ is not uniformly bounded, then as a $\Vir[G']$-module,
$V'\cong V(\a,\b,G'_0,b)$ for some $\a,\b\in \C,\, b\in G'$ and a
subgroup $G'_0$ of $G'$ with $G'=G'_0\oplus \Z b$.  Note that
$\mu_1+\Z g_1\subset \mu_2+G'$. Since $V'_{\mu_1+\sZ g_1}$ is
uniformly bounded, then $g_1\in G'_0$, and hence $V'_{\mu_2+\sZ
g_1}$ is also uniformly bounded by Corollary 2.7.

Thus in both cases,   $V'_{\mu_2+\sZ g_1}$ is a uniformly bounded
$\Vir[\Z g_1]$-sub-quotient of $V_{\mu_2+\sZ g_1}$. By Lemma 3.4,
$V_{\mu_2+\Z g_1}$ is a uniformly bounded $\Vir[\Z g_1]$-module,
forcing $g_1\in G_{\mu_2}$ by Lemma 3.5.

So $G_{\mu_1}\subset G_{\mu_2}$. Symmetrically, we also have
$G_{\mu_2}\subset G_{\mu_1}$. Thus $G_{\mu_1}=G_{\mu_2}, \forall
\mu_1,\mu_2\in \supp V\j\{0\}$.

\medskip
(b)\ Suppose $G^{(0)}\neq G$. We shall prove that $G\cong
G^{(0)}\oplus \Z b $ for some $b\in G$ in three steps.

{\bf Step 1:}
  {\it $G/G^{(0)}$ is torsion-free.}

Otherwise we may choose some $g\in G\j G^{(0)}$ and $k_0\in \N$
such that $k_0g\in G^{(0)}$. Take any $\mu\in \supp V\j\{0\}$.
Since $\rank G>1$ we have a subgroup $A$ of $G$ such that $A\cong
\Z^2$. By considering the nontrivial $\Vir[A]$-module $V_{\mu+A}$
and using Theorem 2.8, we know that there exists  $\l\in
\supp(V)\j\Z g$.

 Then $V_{\l+\sZ g}$ is not
uniformly bounded while $V_{\l+\sZ k_0g}$ is uniformly bounded, by
the definition of $G^{(0)}$. Since $V_{\l+\sZ g}$ is not uniformly
bounded, it has a non-trivial highest or lowest weight irreducible
subquotiens, say non-trivial highest weight  $\Vir[\Z
g]$-sub-quotient $W$. Without loss of generality, we may assume
that $W$ has the highest weight $\lambda\ne0$. Then
$W_{\lambda+\sZ k_0g}$ has a nontrivial highest weigh
sub-quotient. Thus $W_{\lambda+\sZ k_0g}$ is not uniformly
bounded, contradicting the fact that $V_{\l+\sZ k_0g}$ is
uniformly bounded.  Thus $G/G^{(0)}$ is torsion-free. Step 1
follows.

% Then we have that $\dim V_{\l+kg}>\dim V_{\l}$ for some
%  $k\in \Z$ with $\l+kg\neq 0$. Let $I=\{\l, \l+kg\}$, and take
 %$G_I$ as in Lemma 3.2. Set $G'=G_I+\Z g$. Then the $\Vir[G']$
%module $V_{\l+G'}$ has a unique irreducible sub-quotient $V'$ with
%$V'_{\mu}=V_{\mu}, \mu\in I$. $V'$ is not uniformly bounded, then $V'\cong
%V(\a,\b,\G'_0,\b)$ for some $\a,\b\in \C,\, b\in G'$ and a
%subgroup $G'_0$ of $G'$ such that $G'=G'_0\oplus \Z b$.

{\bf Step 2: $\rank(G/G^{(0)})=1$.}

Otherwise we assume that $g_1,g_2\in G\j G^{(0)}$ are independent
modulo $G^{(0)}$, i.e., the subgroup $G'=\langle g_1,g_2\rangle$
is isomorphic to $\Z^2$. Then $G'\cap G^{(0)}=0$. Take any $\l\in
\supp V\j\{0\}$ and consider the $\Vir[G']$-module $V_{\l+G'}$,
which has an irreducible sub-quotient $V'$ with $V'_\l\ne0$.

%Either $V'$ is uniformly bounded or not, we may find some $g\in G'$ such
%that $V'_{\l+\Z g}$ is a uniformly bounded $\Vir[\Z g]$ module as in 1).

If $V'$ is uniformly bounded, then $V'_{\l+\sZ g_1}$ is a
nontrivial uniformly bounded $\Vir[\Z g_1]$-sub-quotient of
$V_{\l+\sZ g_1}$. Thus by Lemma 3.4, $V_{\l+\sZ g_1}$ is a
uniformly bounded $\Vir[\Z g_1]$ module, forcing $g_1\in G^{(0)}$
by Lemma 3.5, a contradiction.

If $V'$ is not uniformly bounded, then $V'\cong V(\a,\b,G'_0,b)$
for some $\a,\b\in \C,\, b\in G'$ and a subgroup $G'_0$ of $G'$
with $G'=G'_0\oplus \Z b$.  Then $V'_{\l+\sZ g}$ is uniformly
bounded for any $g\in G'_0\j\{0\}$, and hence $V_{\l+\sZ g}$ is
also uniformly bounded. Thus $g\in G^{(0)}$, contradicting the
fact that $G'\cap G^{(0)}=\{0\}$. Step 2 follows.

Thus   $G\subset G^{(0)}+\Q g$ for any $g\in G\j G^{(0)}$.

Fix $g\in G\j G^{(0)}$,  we know that $G\subset G^{(0)}+\Q g$.
Since $G/G^{(0)}$ is torsion-free, $G^{(0)}\cap \Q g=\{0\}$. Then
$G=G^{(0)}\oplus G_1$, where $G_1=G\cap \Q g$.

{\bf Step 3:  $G_1 \cong\Z$.}

Otherwise,  $G_1$ would be an infinitely generated abelian group
of rank $1$. Then by Theorem 2.5$'$, the $\Vir[G_1]$ module
$V_{\l+G_1}$ is uniformly bounded for any $\l\in\supp V$, forcing
$G_1\subset G^{(0)}$, contradiction again. Thus we must have
$G_1\cong \Z$. This complete the proof. \qed

\medskip
{\bf Theorem 3.7.}  {\it  Suppose that $V$ is an  irreducible
Harish-Chandra module over the generalized Virasoro algebra
$\Vir[G]$ that is not uniformly bounded. Then $V$ is isomorphic to
$V(\a,\b,G^{(0)},b)$ for some $\a,\b\in \C,\, b\in G$ and a
subgroup $G^{(0)}$ of $G$ with $G=G^{(0)}\oplus \Z b$.}
\par
\medskip
{\it Proof.}\ Note that we have assumed    $\rank G\geq 2$. Since
$V$ is not uniformly bounded, we use Lemma 3.6 to have $G^{(0)}$
and $b$  as described there. Then $C$ trivially acts on $V$.

Take any $\mu\in \supp V\j\{0\}$. Since $\rank G>1$ we have a
subgroup $A$ of $G$ such that $A\cong \Z^2$. By considering the
nontrivial $\Vir[A]$-module $V_{\mu+A}$ and using Theorem 2.8, we
know that there exists  $\l\in \supp(V)\j\Z g$. The $\Vir[\Z
b]$-module $W=V_{\l+\sZ b}$ cannot have any nontrivial uniformly
bounded sub-quotient, for otherwise we would have $b\in G^{(0)}$
by Lemma 3.5, contradicting the definition of $G^{(0)}$. Note that
$0\notin \supp W$. Then by Theorem 2.4, $W=W^+\oplus W^-$, where
$W^+$ is such that $\supp W^+\subset \{\l+kb| k\leq t_0\}$ for
some $t_0\in \Z$ and $W^-$ is such that $\supp W^-\subset \{\l+kb|
k\geq s_0\}$ for some $s_0\in \Z$.

Since $0\notin \l+\Z b$, it is clear that the highest weight
$\l+k_1 b$ of any irreducible highest weight $\Vir[\Z
b]$-sub-quotient of $W$ must satisfy $k_1\le t_0$ and that the
lowest weight $\l+k_2 b$  of any irreducible lowest weight
$\Vir[\Z b]$ sub-quotient of $W$ must satisfy $k_2\ge s_0$.

We want to show that one of $W^+$ and $W^-$ is zero. Otherwise we
may choose $t>t_0$ and $s<s_0$ such that  both $\l+tb$ and $\l+sb$
lie in $\supp W$. Let $I=\{\l+tb,\l+sb\}$, and take $G_I$ as in
Lemma 3.1, and furthermore  we may assume that $\rank G_I>1$. Set
$G'=G_I+\Z b$. The $\Vir[G']$-module $V_{\l+G'}$ has a unique
irreducible sub-quotient $V'$ with $\dim V'_{\mu}=\dim
W_{\mu}=\dim V_{\mu}$ for all $\mu\in I$.

Clearly, $V'$ is not a uniformly bounded $\Vir[G']$-module.
%for otherwise, $V'_{\l+\Z b}$ would be a uniformly bounded $\Vir[\Z b]$
%module, forcing $V_{\l+\Z b}$ being uniformly bounded by Lemma 3.4
%and $b\in G^{(0)}$ by the definition of $G^{(0)}$.
Now by Theorem 2.8, we have that $V'\cong V(\a,\b,G'_0,\epsilon
b)$ for some $\a,\b\in \C,\, \epsilon\in \{1,-1\}$ and a subgroup
$G'_0$ of $G'$ with $G'=G'_0\oplus \Z b$. Thus, $V'_{\l+\sZ b}$ is
either a positively truncated or a negatively truncated $\Vir[\Z
b]$-module, and in both cases $\dim V'_{\l+tb}$ and $\dim
V'_{\l+sb}$ are both nonzero. This implies that either
$W=V_{\l+\sZ b}$ has a highest weight $\Vir[\Z b]$ sub-quotient
with highest weight $\l+kb$ with $k\ge t>t_0$ or has a lowest
weight $\Vir[\Z b]$ sub-quotient with lowest weight $\l+kb$ with
$k\le s<s_0$, contradiction.

Hence, $W=W^+$ or $W=W^-$. Without loss of generality, we assume
that $W=W^+$, that is, $\dim V_{\l+kb}=0, \forall k>t_0$. But
$V_{\l+kb+G^{(0)}}$ is a uniformly bounded $\Vir[G^{(0)}]$-module
for any $k\in \Z$, then we must have that $\dim V_{\l+kb+g}=\dim
V_{\l+kb}=0, \forall k>t_0, g\in G^{(0)}$, provided that
$\l+kb+g\neq 0$.

Let $t_1$ be the largest integer such that $\dim V_{\l+t_1b}\neq
0$.

 If $0\in \l+t_1b+G^{(0)}$, then $\l+t_1b+G^{(0)}=G^{(0)}$.
 Then $\supp(V)\subset -\Z^+b+G^{(0)}$.
Any irreducible $\Vir[G^{(0)}]$-submodule $W$ of $V_{G^{(0)}}$
generates $V$ as $\Vir[G]$-module.  Using PBW theorem we know that
$W=V_{G^{(0)}}$. Since  $V_{G^{(0)}}\supset  V_{\l+t_1b}\neq 0$,
$V_{G^{(0)}}$ is a uniformly bounded nontrivial irreducible
$\Vir[G^{(0)}]$-module. So $V_{G^{(0)}}\cong V'(\a,\b,G^{(0)})$
for some $\a,\b\in\C$. Consequently $V\cong V(\a,\b,G^{(0)},b)$.

If $0\notin \l+t_1b+G^{(0)}$, and $0\notin \l+(t_1+1)b+G^{(0)}$ or
$0\in \l+(t_1+1)b+G^{(0)}$ but $V_0=0$, then
$\Vir_{b+G^{(0)}}V_{G^{(0)}}=0$ where
 $\Vir_{b+G^{(0)}}=\sum_{x\in b+G^{(0)}}\C d_x$. So
$\Vir_{\sN b+G^{(0)}}V_{G^{(0)}}=0$. Thus $\supp(V)\subset
-\Z^+b+G^{(0)}$. Similar discussions yield to the same conclusion
that $V\cong V(\a,\b,G^{(0)},b)$.

If  $0\in \l+(t_1+1)b+G^{(0)}$, and $V_0\ne0$, take nonzero $v\in
V_0$. Then $\Vir_{b+G^{(0)}}v=0$ and $\Vir_{G^{(0)}}v=0$. Thus
$\Vir_{\sZ^+b+G^{(0)}}v=0$. Since $V$ is not trivial, using PBW
theorem and $C=0$ we deduce that $V_{-\sZ b+G^{(0)}}$ is a proper
submodule, contradiction. So this case does not occur.

This proves the theorem.  \qed

\medskip

Combining Theorems 2.2, 2.5$'$,  3.3,  and 3.7 we now have proved
the following classification theorem.

\medskip
{\bf Theorem 3.8.}  {\it  Suppose that $G$ is an arbitrary
additive subgroup of $\C$.}

(a) {\it  If} $\rank G=1$ {\it  and $G\not\cong \Z$, then any
nontrivial irreducible Harish-Chandra module over $\Vir[G]$ is  a
module of the intermediate series.}

(b) {\it If  $G\cong \Z$,  then any nontrivial irreducible
Harish-Chandra module over $\Vir[G]$ is  a module of the
intermediate series, a highest weigh module or a lowest weight
module.}

(c) {\it  If} $\rank G>1$, {\it then any nontrivial irreducible
Harish-Chandra module over $\Vir[G]$ is either a module of the
intermediate series or isomorphic to $V(\a,\b,G^{(0)},b)$ for some
$\a,\b\in \C,\, b\in G$ and a subgroup $G^{(0)}$ of $G$ with
$G=G^{(0)}\oplus \Z b$.}

\vskip 5pt {\bf Acknowledgement.} The authors would like to thank
sincerely Prof.\ Volodymyr Mazorchuk for scrutinizing  the paper,
pointing out several inaccuracies, and making valuable suggestions
to improve the paper.

\par
\vskip .8cm

\cl{\bf{REFERENCE}}
\begin{description}

\item{[BZ]}  Y. Billig and K. Zhao,  Weight modules over
exp-polynomial Lie algebras, {\it J. Pure Appl. Algebra}, Vol.191,
23-42(2004).

\item{[HWZ]} J. Hu, X. Wang, and K. Zhao, Verma modules over
generalized Virasoro algebras $\Vir[G]$, {\it J. Pure Appl.
Algebra}, 177(2003), no.1, 61-69.

\item{[K1]} I. Kaplansky, The Virasoro algebra, {\it Comm. Math.
Phys.}, 86(1982), no.1, 49-54.

\item{[K2]} I. Kaplansky, {\it Infinite abelian groups}, Revised
edition, The University of Michigan Press, Ann Arbor, Mich. 1969.

\item{[KR]} V.Kac and A. Raina, Bombay lectures on highest weight
representations of infinite dimensional Lie algebras, World Sci.,
Singapore, 1987.

\item{[LZ1]} R. Lu and K. Zhao, Classification of irreducible
weight modules over higher rank Virasoro algeras, {\it Adv.
Math.}, in press.

\item{[LZ2]} R. Lu and K. Zhao, Classification of irreducible
weight modules over the twisted Heisenberg-Virasoro algebra,
math.RT/0510194, 2005.

\item{[M]} V. Mazorchuk, Classification of simple Harish-Chandra
modules over $\Q$-Virasoro algebra, {\it Math. Nachr.} 209(2000),
171-177.

\item{[Ma]} O. Mathieu, Classification of Harish-Chandra modules
over the Virasoro algebra, {\it Invent. Math.} 107(1992), 225-234.

\item{[MP]} C. Martin and A. Piard,  Nonbounded indecomposable
admissible modules over the Virasoro algebra, {\it Lett. Math.
Phys.} 23(1991), no. 4, 319--324.

\item{[MZ]} V. Mazorchuk, K. Zhao,   Classification of simple
weight Virasoro modules with a finite-dimensional weight space,
{\it J. Algebra}, in press.

\item{[PZ]} J. Patera and H. Zassenhaus, The higher rank Virasoro
algebras, {\it Comm. Math. Phys.} 136(1991), 1-14.

\item{[S1]} Y.  Su, Simple modules over the high rank Virasoro
algebras, {\it  Comm. Algebra}, 29(2001), no.5, 2067-2080.

\item{[S2]} Y.  Su, Classification of Harish-Chandra modules over
the higher rank Virasoro  algebras, {\it Comm. Math. Phys.}, {\bf
240} (2003), 539-551.

\item{[SZ]} Y. Su and K. Zhao, Generalized Virasoro and
super-Virasoro algebras and modules of intermediate series, {\it
J. Algebra}, 252(2002), no.1, 1-19.

\end{description}

\end{document}